\documentclass[12pt]{article}
\usepackage{amsfonts}
\usepackage{latexsym, amssymb, amsmath, amscd, amsfonts,mathrsfs}

\newtheorem{theorem}{Theorem}[section]
\newtheorem{lemma}[theorem]{Lemma}

\numberwithin{equation}{section}

\def\qed{\hfill \rule{4pt}{7pt}}

\def\pf{\noindent {\it Proof.} }

\begin{document}
\begin{center}
{\large {\bf $q$-Hook Length Formulas for Signed Labeled Forests}}

\vskip 6mm {\small William Y.C. Chen$^1$, Oliver X.Q. Gao$^2$ and Peter L. Guo$^3$\\[%
2mm] Center for Combinatorics, LPMC-TJKLC\\
Nankai University, Tianjin 300071,
P.R. China \\[3mm]
$^1$chen@nankai.edu.cn,  $^2$oliver@cfc.nankai.edu.cn, $^3$lguo@cfc.nankai.edu.cn \\[0pt%
] }
\end{center}

\begin{abstract}
A signed labeled forest is defined as a (plane) forest  labeled by $\{1,2,\ldots, n\}$ along
with minus signs associated to some vertices. Signed labeled forests can be viewed as an extension of signed permutations. We define the inversion number, the flag major index and the R-major index on signed labeled forests. They can be considered as
type $B$ analogues of the indices for labeled forests introduced by
Bj\"{o}rner and Wachs.  The flag major index for signed labeled forests
is based on the flag major index on signed permutations introduced by Adin and Roichman, whereas the R-major index for signed labeled forests is
 based on the  R-major index that we introduce for signed permutations, which is closely related to the major defined by  Reiner. We obtain $q$-hook length formulas by $q$-counting  signed labelings of a given forest with respect to the above indices, from which we see that these three indices are equidistributed for signed labeled forests. Our formulas for the major indices and the inversion number are  type $B$ analogues of the formula due to Bj\"{o}rner and Wachs. We also give a type $D$ analogue  with respect to the  inversion number of  even-signed labeled forests.
\end{abstract}



\section{Introduction}

The inversion number and major index for the permutation  group are two of the most important statistics
which have received remarkable attention  in the combinatorial
literature, see, e.g., \cite{Foata,FandS,Garsia,GandG,Goulden,MacMahon}. Extensions of these two statistics have been intensely studied. The present paper concerns two directions of such extensions: One is toward labeled forests duo to Bj\"{o}rner and Wachs \cite{BW}, and the other is toward the Coxeter groups (mainly types $B$ and $D$), see, e.g., \cite{ABR,ABR1,AR,Biagioli,Biagioli2,BZ,Chow,Reiner}.

For integers $m,n$ $(m\leq n)$, we use $[m,n]$ to denote the interval $\{m,m+1,\ldots,n\}$.  Denote by $S_n$ the permutation group on $[1,n]$. Note that we use $[1,n]$, instead of $[n]$, to denote the set $\{1, 2, \ldots n\}$
before we shall use $[n]$ to denote the $q$-number $1+q+q^2+\cdots + q^{n-1}$.
We shall
represent a permutation $\pi\in S_n$ in one-line notation
$\pi=\pi_1\cdots\pi_n$. Then the descent set of $\pi$ is defined by
\begin{equation*}
\mathrm{Des}(\pi)=\{i\colon\ 1\leq i\leq n-1,\, \pi_i>\pi_{i+1}\}.
\end{equation*}
The inversion number and major index of $\pi$ are defined by
\begin{align*}
\mathrm{inv}(\pi)&=|\{(i,j)\colon\ 1\leq i<j\leq n,\, \pi_i>\pi_j\}|,\\[5pt]
\mathrm{maj}(\pi)&=\sum_{i\in \mathrm{Des}(\pi)}i.
\end{align*}
The following classical formula is  duo to  MacMahon \cite{MacMahon}
\begin{equation}\label{MacMahon1}
\sum_{\pi \in S_n}q^{\mathrm{inv}(\pi)}=[n]!=\sum_{\pi \in
S_n}q^{\mathrm{maj}(\pi)},
\end{equation} where $[n]=1+q+\cdots+q^{n-1}$, and $[n]!=[1][2]\cdots
[n]$.

The inversion number and the major index on permutations have been generalized
by Bj\"{o}rner and Wachs \cite{BW} to labeled forests.
Let $F$ be a (plane) forest with vertex set $V(F)$. The reason that we consider plane forests is that every vertex can be viewed as having
a unique position in the
 sense that all the vertices are implicitly labeled.
 A labeling $w$ of $F$ is a
bijection
\[w\colon\ V(F)\longrightarrow [1,n].\] For each vertex $u\in F$, the
hook length of $u$, denoted by $h_u$, is the size of
the subtree rooted at $u$. When $F$ is considered  as a
poset with roots at the top, the hook length of $u$ equals the  cardinality of
the principle ideal $\{v\in F\colon\ v\leq_F u\}$ where $\leq_F$ is the order relation.
Bj\"{o}rner and Wachs \cite{BW} defined the descent set of a labeled forest
as given below,
\[
\mathrm{Des}(F,w)=\{u\in F\colon\ w(u)>w(v),\, \mbox{$v$ is the
parent of $u$}\},
\]
and
\begin{align}
\mathrm{inv}(F,w)&=|\{(u,v)\colon\ u>_Fv,\, w(u)<w(v)\}|,\label{BWinv}\\[5pt]
\mathrm{maj}(F,w)&=\sum_{u\in \mathrm{Des}(F,w)}h_u.\label{BWmaj}
\end{align}
If $F$ is a linear tree, then we get a  permutation by reading
 the labels bottom up. The descent set, the inversion number and the major index for
a linear tree  coincide with the corresponding indices
for permutations. Bj\"{o}rner and Wachs \cite{BW} derived the
following $q$-hook length formula by $q$-counting all labelings of a
fixed forest with respect to the inversion number   and the major index, which reduces to the formula (\ref{MacMahon1}) by
restricting $F$ to be a linear tree.

\begin{theorem}[\mdseries{Bj\"{o}rner and Wachs \cite{BW}, Theorem
1.3}]\label{TBJ1} Let $F$ be a forest of size $n$. Then
\begin{equation}\label{BW}
\sum_{(F,w)}q^{\mathrm{inv}(F,w)}=\frac{n!}{\prod_{u\in
F}h_u}\prod_{u\in
F}[h_u]=\sum_{(F,w)}q^{\mathrm{maj}(F,w)},\end{equation} where
$w$ ranges over all labelings of $F$.
\end{theorem}

In this paper, we shall be concerned with signed labelings of a (plane)
forest.
The signed permutation group $B_n$ is  the group of  bijections $\sigma$ on the set $[-n, n]\backslash \{0\}$ such that
\[\sigma(-i)=-\sigma(i)\]
for $i\in [-n, n]\backslash \{0\}$. Recall that $B_n$ is also known as the hyperoctahedral group of rank $n$, or the Coxeter group of type $B_n$. For $\sigma\in B_n$, we write $\sigma$ in the one-line notation $\sigma=\sigma_1\cdots \sigma_n$, where $\sigma_i=\sigma(i)$ for $i\in [1,n]$. In the language of Coxeter groups, $B_n$ is the Coxeter group of type $B_n$ with respect to the generating set
$\{s_0,s_1,\ldots,s_{n-1}\}$, where $s_i$,
$i=1,2,\ldots,n-1$, are the simple transpositions
\[s_i=(1,\dots,i-1,i+1,i,i+2,\ldots,n)\]and $s_0$ is the sign
change\[s_0=(-1,2,\ldots,n).\]

The length function of an element $\pi$ in a Coxeter group, denoted by $\ell(\pi)$, is the minimum number of
generators that occur in its factorizations, see, Bj\"{o}rner and Brenti \cite{BB}, namely,
\[\ell(\pi)=\min\{r\colon\ \pi=s_{i_1}s_{i_2}\cdots s_{i_r},\
\mbox{$s_{i_1},s_{i_2},\ldots, s_{i_r}$ are generarors}\}.\]

For $B_n$, there exists a simple combinatorial interpretation for the length function.  Let $\sigma=\sigma_1\cdots\sigma_n$ be a signed permutation, and define
\begin{align}
\mathrm{n_1}(\sigma)&=|\{\sigma_i\colon\ 1\leq i\leq n,\ \sigma_i<0\}|,\\[5pt]
\mathrm{n_2}(\sigma)&=\left|\big\{\{i,j\}\colon\ 1\leq i,j\leq n,\
\sigma_i+\sigma_j<0\big\}\right|.\label{N2}
\end{align}
Then the length function of $\sigma$ is given by
\begin{equation}\label{LFB}
\ell_B(\sigma)=\mathrm{inv}(\sigma)+\mathrm{n_1}(\sigma)+\mathrm{n_2}(\sigma),
\end{equation}which can be seen as the inversion number for signed permutations, see, Biagioli \cite{Biagioli}, or Bj\"{o}rner and Brenti \cite{BB}.
The following length generating function is well-known.

\begin{theorem}[\mdseries{Humphreys \cite[Section 3.15]{Humphreys}}]\label{LGFB}
\begin{equation}\label{lenB}
\sum_{\sigma\in B_n}q^{\ell_B(\sigma)}=[2][4]\cdots[2n].
\end{equation}
\end{theorem}

As a subgroup of $B_n$, the group of even-signed permutations, i.e., permutations with an even number of minus signs, is  denoted by $D_n$. It is
well known that  $D_n$ is the Coxeter group  respect to the generating set $\{t_0, s_1,\ldots, s_{n-1}\}$, where $s_i$ for $1\leq i\leq n-1$ are defined as above and
\[t_0=(-2,-1,3,\ldots,n).\]
For any signed permutation $\sigma$, let
\[\mathrm{Neg}(\sigma)=\{1\leq i\leq n\colon \sigma_i<0\}.\]
The length function for $D_n$ can be computed by the following combinatorial
formula
\begin{equation}\label{lengthD}
\ell_D(\sigma)=\mathrm{inv}(\sigma)-\sum_{i\in \mathrm{Neg}(\sigma)}\sigma_i-\mathrm{n}_1(\sigma),
\end{equation}which can be considered as the inversion number of an even-signed permutation, see, Bj\"{o}rner and Brenti \cite[Section 8.2]{BB}, or Biagioli \cite{Biagioli2}.
It can be checked that
\[-\sum_{i\in \mathrm{Neg}(\sigma)}\sigma_i=\mathrm{n}_1(\sigma)+\mathrm{n}_2(\sigma).\]
So (\ref{lengthD}) can be reformulated as
\begin{equation}\label{LFD}
\ell_D(\sigma)=\mathrm{inv}(\sigma)+\mathrm{n}_2(\sigma).
\end{equation}

The length generating function for $D_n$ is given by the following formula.

\begin{theorem}[\mdseries{Humphreys \cite[Section 3.15]{Humphreys}}]\label{LGFD}
\begin{equation}\label{lenD1}
\sum_{\sigma\in D_n}q^{\ell_D(\sigma)}=[2][4]\cdots[2n-2][n].
\end{equation}
\end{theorem}

Statistics on Coxeter groups  that are equiditributed with  the length function
are called Mahonian.
 An important statistic on $B_n$ is the flag major index introduced by Adin and Roichman \cite{AR},
which is defined in terms of Coxeter elements and can be expressed combinatorially as
\begin{equation}\label{flagmaj}
\mathrm{fmaj}(\sigma)=2\,\mathrm{maj}(\sigma)+\mathrm{n_1}(\sigma).
\end{equation}

The second Mahonian major statistic on $B_n$ is the negative major index introduced by  Adin,  Brenti and  Roichman \cite{ABR}, which has the following combinatorial description
\begin{equation}\label{nmaj}
\mathrm{nmaj}(\sigma)=\mathrm{maj}(\sigma)+\mathrm{n}_1(\sigma)+\mathrm{n}_2(\sigma).
\end{equation}

There is another Mahonian  statistic which is based on  the major index  defined by Reiner \cite{Reiner}. Under the following order
\begin{equation}\label{order}
1<\cdots<n<-n<\cdots<-1,
\end{equation}
the descent set $\mathrm{Des}_R(\sigma)$ is defined as
\begin{equation*}
\mathrm{Des}_R(\sigma)=\{i\in [1,n]\colon\, \sigma_i>\sigma_{i+1}\},
\end{equation*}
under the assumption that $\sigma_{n+1}=n$.
Then the major index $\mathrm{maj}_R(\sigma)$ is given by
\begin{equation}\label{Rmaj}
\mathrm{maj}_R(\sigma)=\sum_{i\in \mathrm{Des}_R(\sigma)}i.
\end{equation}
Reiner \cite{Reiner} has shown that
\begin{equation}\label{rg}
\sum_{\sigma\in B_n}t^{\mathrm{n_1}(\sigma)}q^{\mathrm{maj}_R(\sigma)}=(1+tq)^n[n]!.
\end{equation}
While the index $\mathrm{maj}_R$ is not Mahonian,
as observed by Biaginoli and Zeng \cite{BZ}, from   (\ref{rg}) it follows that
the index
\begin{equation}\label{zeng}
2\,\mathrm{maj}_R(\sigma)-\mathrm{n}_1(\sigma)
\end{equation} is  equidistributed with the flag major index (\ref{flagmaj}).
They also gave a proof of this fact by justifying the following relation
\begin{equation*}
\mathrm{maj}_R(\sigma)=\mathrm{maj}(\sigma)+\mathrm{n}_1(\sigma).\end{equation*}

We next define a new Mahonian index called R-major index for signed permutations which relies on the major index $\mathrm{maj}_B$ with respect to
 the natural order
 \begin{equation}\label{norder}
  -n < \cdots <-1<0< 1<\cdots < n.
  \end{equation}
The index $\mathrm{maj}_B$  can be shown to be isomorphic to Reiner's major index (\ref{Rmaj}) defined with respect to the order (\ref{order}).
The descent set $\mathrm{Des}_B(\sigma)$ is described in terms of the natural order; that is,
\begin{equation}
\mathrm{Des}_B(\sigma)=\{i\in [1,n]\colon\, \sigma_i>\sigma_{i+1}\},
\end{equation} where $\sigma_{n+1}=0$.  Let
\begin{equation}\label{PRmaj}
\mathrm{maj}_B(\sigma)=\sum_{i\in \mathrm{Des}_B(\sigma)}i,
\end{equation}
and let $\mathrm{p}(\sigma)$ be the number of positive entries of $\sigma$. Then the R-major index is given by
\begin{equation}\label{positive}
\mathrm{rmaj}(\sigma)=2\,\mathrm{maj}_B(\sigma)-\mathrm{p}(\sigma).
\end{equation}

It can be shown that the  index $\mathrm{maj}_B(\sigma)$ in (\ref{PRmaj}) is isomorphic to the index $\mathrm{maj}_R(\sigma)$ in (\ref{Rmaj}) of Reiner, see Lemma \ref{lastlemma}. Under this correspondence, the index $\mathrm{rmaj}(\sigma)$ in (\ref{positive}) is isomorphic to the index $2\,\mathrm{maj}_R(\sigma)-\mathrm{n}_1(\sigma)$ in (\ref{zeng}).
It should be mentioned that we prefer the  R-major index formulation
 based on the the natural order (\ref{norder})
because it is consistent with the order of the flag major index and the negative
major index and it seems to be easier to describe.

In this paper, we extend the three indices on signed permutations, i.e., the inversion number, the flag major index,
and the R-major index,  to signed labeled forests. The inversion number for signed labeled forests is motivated by the inversion number (\ref{BWinv}) and the length function (\ref{LFB}).  The flag major index for
signed labeled forests is an extension of the major indices (\ref{BWmaj}) for forests and (\ref{flagmaj}) for signed permutations, whereas the R-major index for signed labeled forests stems from the major indices  (\ref{BWmaj}) and (\ref{positive}).

 We obtain $q$-hook length formulas by $q$-counting  signed labelings of a given forest with respect to the above indices, from which we see that these three indices are equidistributed for signed labeled forests. Our formulas for the major indices and the inversion number are  type $B$ analogues of the formula due to Bj\"{o}rner and Wachs. We also define the inversion number on even-signed labeled forests in connection with the inversion number (\ref{BWinv}) and the length function (\ref{LFD}) for even-signed permutations.

This paper is organized as follows. In Section 2, we give the definitions of the inversion number and the two major indices on singed labeled forests. The main results are described in this section. We also include a sketch of the proof of the $q$-hook length  formula for the flag major index for signed labeled forests.  Section 3 is   devoted to the
proofs of the generating function formulas for the inversion numbers on signed labeled forests and even-signed labeled forests. In Section 4, we consider the generating function of the R-major index. To this end, we define
$(P,w)$-partitions of type $B$ which can be viewed as a type $B$ extension of the usual
$(P,w)$-partitions as introduced by Stanley \cite{Stanley2}. In Section
5, we give a bijection which establishes the connection  between
the flag  major index and the R-major index on signed labeled forests. Section 6 gives some further questions.

\section{Main results}

In this section, we give the definitions of the inversion number, the flag
 major index and the R-major index on singed (even-signed) labeled forests as aforementioned  in the introduction, and outline the main results of this paper.

Let $F$ be a forest.
Denote by $B_n(F)$ (resp., $D_n(F)$) the set of  signed (resp., even-signed)
labeled forests with the underlying forest $F$. We use $(F,w)$ to
denote a signed labeling of $F$.  The motivation to consider types $B$ and $D$ analogues for signed  labeled
forests  is  the observation that the number $\mathrm{n_2}(\sigma)$ in (\ref{N2}) has a
natural extension to  signed labeled forests; that is, for
$(F,w)\in B_n(F)$, we may define
\begin{equation}\label{Bextension}
\mathrm{n}_2(F,w)=|\{(u,v)\colon u<_Fv, w(u)+w(v)<0\}|.
\end{equation}
Let $\mathrm{n}_1(F,w)$ be the number of negative labels of $w$. So we define the inversion number
for $(F,w)\in B_n(F)$ as
\begin{equation}\label{SLFinv}
\mathrm{inv}_B(F,w)=\mathrm{inv}(F,w)+\mathrm{n}_1(F,w)+\mathrm{n}_2(F,w),
\end{equation}
and  for $(F,w)\in D_n(F)$ the inversion number is defined by
\begin{equation}\label{ESLFinv}
\mathrm{inv}_D(F,w)=\mathrm{inv}(F,w)+\mathrm{n}_2(F,w).
\end{equation}

The first two theorems assert that the inversion numbers (\ref{SLFinv}) and (\ref{ESLFinv}) lead to types $B$ and $D$ analogues of the length generating functions (\ref{lenB}) and (\ref{lenD1}) respectively.

\begin{theorem}\label{ThminvB}
Let $F$ be a forest of size $n$. Then
\begin{equation}\label{lenB2}
\sum_{(F,\,w)\in\, B_n(F)}q^{\mathrm{inv}_B(F,\,w)}=\frac{n!}{\prod_{u\in F}h_u}\prod_{u\in
F}[2h_u].\end{equation}
\end{theorem}

\begin{theorem}\label{ThminvD}
Let $F$ be a forest of size $n$. Then
\begin{equation}\label{sinvD}
\sum_{(F,\,w)\in\, D_n(F)}q^{\mathrm{inv}_D(F,\,w)}=\frac{n!}{2\prod_{u\in
F}h_u}\prod_{u\in F}(1+q^{h_u-1})[h_u].\end{equation}
\end{theorem}

Note that Theorem \ref{LGFB} and Theorem \ref{LGFD} can be deduced respectively from Theorem \ref{ThminvB} and Theorem \ref{ThminvD} by restricting  $F$ to  a linear order.

We next turn to two major  indices for signed labeled forests. The first  is based on the flag major index for signed permutations, namely, for
$(F,w)\in B_n(F)$, we define
\begin{equation}\label{SLFfmaj}
\mathrm{fmaj}(F,w)=2\,\mathrm{maj}(F,w)+\mathrm{n}_1(F,w).
\end{equation}
The following theorem shows that the flag major index (\ref{SLFfmaj}) is equidistributed with the inversion number (\ref{SLFinv}) for signed labeled
forests.

\begin{theorem}\label{Tfmaj}
Let $F$ be a forest of size $n$. Then
\begin{equation}\label{fmaj}
\sum_{(F,\,w)\in\, B_n(F)}q^{\mathrm{fmaj}(F,\,w)}=\frac{n!}{\prod_{u\in
F}h_u}\prod_{u\in F}[2h_u].
\end{equation}
\end{theorem}

On the other hand, it will be seen  that Theorem \ref{Tfmaj} can be deduced from  Theorem \ref{TBJ1} and
a decomposition of the signed permutation group. We sketch the proof below.
Let \[T_n=\{\sigma_1\cdots\sigma_n\in B_n\colon\, \sigma_1<\cdots<\sigma_n\}.\]
Then the signed permutation group $B_n$ has the following decomposition
\begin{equation}\label{decomposition}
B_n=\bigcup_{\pi\in S_n}\{\sigma\pi\colon\, \sigma\in T_n\},
\end{equation}
which is the multiplicative
decomposition of a Coxeter group into a parabolic subgroup and its minimal coset
representatives, see, e.g., Bj\"{o}rner and Brenti \cite{BB} or Humphreys \cite{Humphreys}. Such a decomposition has been used by Adin et al. \cite[Proposition 3.1]{ABR} to prove the Mahonian property of the negative major index (\ref{nmaj}).

\noindent\textit{Sketch of the proof of Theorem \ref{Tfmaj}.} Let $w_0$ be
 a decreasing labeling of $F$, i.e., a labeling of $F$ such that $w_0(u)>w_0(v)$ whenever $u>v$. We now define a bijection  $\psi\colon (F,w)\longmapsto \sigma$  from  $B_n(F)$ to $B_n$ such that $\sigma_i=w(i)$ for $1\leq i\leq n$. Let
\[T_n(F)=\{\psi^{-1}(\sigma)\colon\, \sigma\in T_n\},\]
and let $U_n(F)=\psi^{-1}(S_n)$ be the set of ordinary  labelings of $F$. From the decomposition (\ref{decomposition}), we get
\begin{equation}\label{decompsition2}
B_n(F)=\bigcup_{w\in U_n(F)}\{\tau w\colon\, \tau\in T_n(F)\},
\end{equation}where $\tau w(i)=\tau (w(i))$ for any vertex $i\in F$.

Let $u$ be the parent of $v$, and let  $w_0(u)=i_{u}$ and $w_0(v)=i_v$. It is not hard to verify that $\tau w(i_u)>\tau w(i_v)$ if and only if $w(i_u)>w(i_v)$. This  implies that $\tau w$ and $w$ have the same descent set. Similarly, we see that  $\mathrm{n_1}(F,\,\tau w)=\mathrm{n_1}(F,\,\tau)$. Thus
\begin{equation}\label{key}
\begin{split}
\mathrm{fmaj}(F,\,\tau w)&=2\,\mathrm{maj}(F,\,\tau w)+\mathrm{n_1}(F,\,\tau w)\\
&=2\,\mathrm{maj}(F,\,w)+\mathrm{n_1}(F,\,\tau).
\end{split}
\end{equation}
Again, by the decomposition (\ref{key}) and  Theorem \ref{TBJ1}, we obtain
\begin{align*}
\sum_{(F,\,w)\in B_n(F)}q^{\mathrm{fmaj}(F,\,w)}
&=\sum_{w\in U_n(F)}\sum_{\tau\in T_n(F)}q^{2\,\mathrm{maj}(F,\,w)+\mathrm{n_1}(F,\,\tau)} \\[5pt]
&=\sum_{\tau\in T_n(F)}q^{\mathrm{n_1}(F,\,w)}\ \frac{n!}{\prod_{u\in
F}h_u}\prod_{u\in
F}[h_u]_{q^2} \\[5pt]
&=(1+q)^n\frac{n!}{\prod_{u\in
F}h_u}\prod_{u\in
F}[h_u]_{q^2}\\[5pt]
&=\frac{n!}{\prod_{u\in
F}h_u}\prod_{u\in
F}(\left[h_u]_{q^2}(1+q)\right)\\[5pt]
&=\frac{n!}{\prod_{u\in
F}h_u}\prod_{u\in
F}[2h_u],
\end{align*}
as desired. \qed

The second major index  for  signed labeled forests is the R-major index. For $(F,w)\in B_n(F)$,  we define the
descent set  as
\[\mathrm{Des}_B(F,w)=\mathrm{Des}(F,w)\bigcup\,\{u\in F\colon\, \mbox{$u$ is a root of $F$ with a positive label}\}.\]
Let
\begin{equation}\label{major2}
\mathrm{maj}_B(F,w)=\sum_{u\in \mathrm{Des}_B(F,w)}h_u.
\end{equation}
Denote by $\mathrm{p}(F,w)$ the number of positive labels of $w$. Then the R-major index is defined by
\begin{equation}\label{SLFrmaj}
\mathrm{rmaj}(F,w)=2\,\mathrm{maj}_B(F,w)-\mathrm{p}(F,w).
\end{equation}
The following theorem shows that the R-major index is
equidistributed with the inversion number  for signed labeled forests.

\begin{theorem}\label{Tnmaj}
Let $F$ be a forest of size $n$. Then
\begin{equation}\label{nnmaj}
\sum_{(F,\,w)\in\, B_n(F)}q^{\mathrm{rmaj}(F,\,w)}=\frac{n!}{\prod_{u\in
F}h_u}\prod_{u\in F}[2h_u].
\end{equation}
\end{theorem}

\section{The inversion numbers for types $B$ and $D$}

In this section, we are aimed to prove Theorems \ref{ThminvB} and \ref{ThminvD}, both of which  can be deduced from  the following  theorem involving the weighted counting  of signed labeled forest with respect to the number of inversions and the number of negative labels.

\begin{theorem}\label{thmsinv}
Let $F$ be a forest of size $n$. Then
\begin{equation}\label{sinv}
\sum_{(F,\,w)\in\, B_n(F)}t^{\mathrm{n}_1(F,\,w)}q^{\mathrm{inv}_B(F,\,w)}=\frac{n!}{\prod_{u\in
F}h_u}\prod_{u\in F}(1+tq^{h_u})[h_u].\end{equation}
\end{theorem}

It is clear that  the constant terms on both sides of (\ref{sinv}) lead to
the $q$-hook length formula of Bj\"orner and Wachs as stated in Theorem \ref{TBJ1}.

\pf We proceed by induction on $n$. Assume that  $F$ is a forest consisting of $k$ trees $T_1,\ldots,T_k$. Let us consider the following two cases.

Case 1. Suppose that $k>1$ and  $T_i$ has $n_i$ vertices for  $1\leq i\leq
k$. Each signed labeling of $F$ corresponds  to a
$k$-tuple consisting of a distribution of the $n$ labels (some of
which may have minus signs) into $k$ trees with each tree $T_i$
receiving $n_i$ labels. We see that
\begin{equation*}
\sum_{(F,\,w)}t^{\mathrm{n}_1(F,w)}q^{\mathrm{inv}_B(F,w)}={n\choose
{n_1,\ldots,n_k}}\prod_{i=1}^k
\sum_{(T_i,\,w_i)}t^{\mathrm{n}_1(T_i,w_i)}q^{\mathrm{inv}_B(T_i,w_i)}.
\end{equation*}
Applying the induction hypothesis yields
\begin{align*}
\sum_{(F,\,w)}t^{\mathrm{n}_1(F,w)}q^{\mathrm{inv}_B(F,w)}&={n\choose
{n_1,\ldots,n_k}}\prod_{i=1}^k\frac{n_i!}{\prod_{u\in
T_i}h_u}\prod_{u\in T_i}(1+tq^{h_u})[h_u]\\
&=\frac{n!}{\prod_{u\in F}h_u}\prod_{u\in F}(1+tq^{h_u})[h_u].
\end{align*}

 Case 2. Suppose that $k=1$. Let $F'$ be the
forest obtained from $F$ by removing the root $u$. Every signed
labeling $w$ of $F$ corresponds  to a pair $(w(u),w')$, where $w'$ is the signed labeling of $F'$ induced by $w$. Since $w$ is a signed labeling, we need to
consider the two subcases $w(u)>0$ and $w(u)<0$.

Assume that $w(u)>0$. Then let $w(u)=i$, where $1\leq i\leq n$. We get
\begin{align*}
\mathrm{inv}(F,w)&=\mathrm{inv}(F',w')+|\{v\in F'\colon w'(v)>i\}|,\\[5pt]
\mathrm{n}_1(F,w)&=\mathrm{n}_1(F',w'),\\[5pt]
\mathrm{n}_2(F,w)&=\mathrm{n}_2(F',w')+|\{v\in F'\colon w'(v)+i<0\}|.
\end{align*}
It is easy to check
\[\{v\in F'\colon w'(v)>i\}\cup\{v\in F'\colon w'(v)+i<0\}=\{v\colon i<|w'(v)|\leq n\},\] where $\cup$ denotes the disjoint union. Thus we find
\[\mathrm{inv}_B(F,w)=\mathrm{inv}_B(F',w')+n-i.\]

 Assume that $w(u)<0$. Then let $w(u)=-i$, where $1\leq i\leq n$. We have
\begin{align*}
\mathrm{inv}(F,w)&=\mathrm{inv}(F',w')+|\{v\in F'\colon w'(v)>-i\}|,\\[5pt]
\mathrm{n}_1(F,w)&=\mathrm{n}_1(F',w')+1,\\[5pt]
\mathrm{n}_2(F,w)&=\mathrm{n}_2(F',w')+|\{v\in F'\colon w'(v)-i<0\}|.
\end{align*}
We need to determine the numbers $c_1$ and $c_2$ defined by
\[ c_1=|\{v\in F'\colon w'(v)>-i\}|,\quad   c_2=|\{v\in
F'\colon w'(v)-i<0\}|.\]
\begin{itemize}
\item[(i).] Consider the index  $j$ with $i<j \leq n$. If $j$ is in the labeling set of $w$,  then $j$
contributes $1$ to $c_1$ and $0$ to  $c_2$. If $-j$
is in the labeling set of $w$, then $-j$ contributes $0$ to
$c_1$ and $1$ to  $c_0$.
\item[(ii).]
Consider the index $j$ with $1\leq j<i$. If $j$ is in the labeling set of $w$, then $j$
contributes $1$  to both $c_1$ and $c_2$.  If $-j$ is in the labeling set of $w$, then $-j$ also
contributes $1$ to both $c_1$ and $c_2$.
\end{itemize}
\allowdisplaybreaks
It follows that
\begin{equation*}
|\{v\in F'\colon w'(v)>-i\}+|\{v\in F'\colon
w'(v)-i<0\}|=(n-i)+2(i-1)=n-2+i,
\end{equation*}
which gives  \begin{align*} \mathrm{inv}_B(F,w)&=\mathrm{inv}_B(F',w')+c_1+c_2+1\\[5pt]
&=\mathrm{inv}_B(F',w')+n-1+i.\end{align*} Combining the above two cases $w(u)>0$ and $w(u)<0$ and using the induction hypothesis, we deduce that
\begin{align*}
&\sum_{(F,\,w)}t^{\mathrm{n}_1(F,\,w)}q^{\mathrm{inv}_B(F,\,w)}\\[5pt]
=&\sum_{i=1}^nq^{n-i}\sum_{(F',\,w')}t^{\mathrm{n}_1(F',\,w')}q^{\mathrm{inv}_B(F',\,w')}
+\sum_{i=1}^ntq^{n-1+i}\sum_{(F',\,w')}t^{\mathrm{n}_1(F',\,w')}q^{\mathrm{inv}_B(F',\,w')}\\[5pt]
=&(1+q+\cdots+q^{n-1}+t(q^n+\cdots+q^{2n-1}))\sum_{(F',\,w')}t^{\mathrm{n}_1(F',\,w')}q^{\mathrm{inv}_B(F',\,w')},\\[5pt]
=&(1+tq^n)[n]\sum_{(F',\,w')}t^{\mathrm{n}_1(F',\,w')}q^{\mathrm{inv}_B(F',\,w')}\\[5pt]
=&(1+tq^n)[n]\frac{(n-1)!}{\prod_{u\in F'}h_u}\prod_{u\in
F'}(1+tq^{h_u})[h_u]\\[5pt]
=&\frac{n!}{\prod_{u\in F}h_u}\prod_{u\in F}(1+tq^{h_u})[h_u],
\end{align*}
as claimed.  \qed

Once Theorem \ref{thmsinv} is established, it is easy to derive Theorem \ref{ThminvB} by taking $t=1$ in (\ref{sinv})
and by observing that
\[(1+q^n)[n]=[2n].\]

We are now ready to give a proof of Theorem \ref{ThminvD}.

\noindent\textit{Proof of Theorem \ref{ThminvD}.}  Let
\begin{equation*}
D_n(t,q)=\sum_{(F,\,w)\in B_n(F)}t^{\mathrm{n}_1(F,\,w)}q^{\mathrm{inv}(F,w)+\mathrm{n}_2(F,w)}.
\end{equation*}
By Theorem \ref{thmsinv} we obtain
\begin{align*}
D_n(t,q)=&\sum_{(F,\,w)\in B_n(F)}(t/q)^{\mathrm{n}_1(F,\,w)}q^{\mathrm{inv}(F,w)+\mathrm{n}_1(F,w)+\mathrm{n}_2(F,w)}\\[5pt]
=&\sum_{(F,\,w)\in B_n(F)}(t/q)^{\mathrm{n}_1(F,\,w)}q^{\mathrm{inv}_B(F,\,w)}\\[5pt]
=&\frac{n!}{\prod_{u\in F}h_u}\prod_{u\in F}(1+tq^{h_u-1})[h_u].
\end{align*}
Since $F$ has at least one leaf, there must exist a factor $(1+tq^{h_u-1})=(1+t)$  in $D_n(t,q)$. Hence $D_n(-1,q)=0$, or, equivalently,
\begin{equation}\label{even}
\sum_{\mbox{$i$ is even}}[t^i]D_n(t,q)=\sum_{\mbox{$i$ is odd}}[t^i]D_n(t,q),\end{equation}where $[t^i]$ denotes the coefficient of $t^i$ in $D_n(t,q)$.
Thus
\begin{align*}
\sum_{(F,\,w)\in\, D_n(F)}q^{\mathrm{inv}_D(F,\,w)}&=\sum_{\mbox{$i$ is even}}[t^i]D_n(t,q)\\[5pt]
&=\frac{D_n(1,q)}{2},\end{align*} which coincides with the formula (\ref{sinvD}). \qed

It would be interesting to give a combinatorial interpretation for (\ref{even}).

\section{The R-major index for signed labeled forests}

In this section, we shall present a proof of Theorem \ref{Tnmaj}. This theorem will be deduced from the following more general formula.

\begin{theorem}\label{Thm1}
Let $F$ be a forest of size  $n$. Then
\begin{equation}\label{main}
\sum_{(F,w)\in\, B_n(F)}t^{\mathrm{p}(F,w)}q^{\mathrm{maj}_B(F,w)}=\frac{n!}{\prod_{u\in
F}h_u}(1+tq)^n \prod_{u\in F}[h_u].\end{equation}
\end{theorem}

 In fact, the proof of the above Theorem relies on
 the formula for $q$-counting the major index $\mathrm{maj}_B$ of  linear extensions of a signed labeled forest, which can be considered as type $B$
  analogue of the $q$-hook length formula  of Bj\"{o}rner and Wachs \cite{BW} for the $q$-counting of the major index of linear extensions.

Let us recall the definition of linear extensions of poset. For a poset $P$ with $n$ elements $x_1,\ldots,x_n$, linear extensions of $P$ can be seen as  permutations $x_{i_1}\cdots x_{i_n}$ such that   $x_{i_j}<_Px_{i_k}$ implies $j<k$.  A labeling of $P$ is a bijection from $\{x_1,\ldots, x_n\}$ to $[1,n]$. For a poset  $P$ with a labeling $w$, we usually use a permutation $w(x_{i_1})\cdots w(x_{i_n})$  to represent  linear extension $x_{i_1}\cdots x_{i_n}$. Denote by $\mathcal{L}(P,w)$ the set of all such permutations. Bj\"{o}rner and Wachs \cite{BW} obtained  the following generating function for the major index of linear extensions of any given labeled forest, which reduces to the result of Stanley \cite{Stanley2} when the labeling is decreasing.

\begin{theorem}[\mdseries{Bj\"{o}rner and Wachs \cite{BW}, Theorem
1.2}]\label{Thm3}Let $F$ be a forest of size $n$, and $w$ a labeling
of $F$. Then
\begin{equation}\sum_{\pi \in
\mathcal{L}(F,w)}q^{\mathrm{maj}(\pi)}=q^{\mathrm{maj}(F,w)}\frac{[n]!}{\prod_{u_\in
F}[h_u]}. \end{equation}
\end{theorem}

When $w$ is a signed labeling,  $\mathcal{L}(F,w)$, as defined above, is a set of signed permutations.
We obtain the following theorem which can be viewed as a type $B$ analogue of Theorem \ref{Thm3}.

\begin{theorem}\label{Le1}
Let $F$ be a forest of size $n$, and $w$ a signed labeling of $F$. Then
\begin{equation}\label{Le1equ}
\sum_{\sigma \in
\mathcal{L}(F,w)}q^{\mathrm{maj}_B(\sigma)}=q^{\mathrm{maj}_B(F,w)}\frac{[n]!}{\prod_{u\in
F}[h_u]}.\end{equation}
\end{theorem}

Though Theorem \ref{Le1} is a type $B$ analogue of Theorem \ref{Thm3}, its proof does not seem to be straightforward.  To prove Theorem \ref{Thm3}, Bj\"{o}rner and Wachs \cite{BW} defined the recursive labelings for forests. They first derived  the following $q$-hook length formula
\[\sum_{\pi \in
\mathcal{L}(F,w)}q^{\mathrm{inv}(\pi)}=q^{\mathrm{inv}(F,w)}\frac{[n]!}{\prod_{u_\in
F}[h_u]},\]
where $w$ is a recursive labeling, see Bj\"{o}rner and Wachs \cite[Theorem 1.1]{BW}. Then they proved that Foata's bijection (see, Foata \cite{Foata}) is invariant when restricted to the set of linear extensions of a forest with a recursive labeling, see Bj\"{o}rner and Wachs \cite[Theorem 2.2]{BW}. Moreover they observed that the inversion number and the major index of a labeled forest are equal for a recursive labeling, see Bj\"{o}rner and Wachs \cite[Lemma 2.3]{BW}. For recursive labelings,  they obtained Theorem \ref{Thm3}. Finally,
 they used a equivalence relation on labelings to extend the special case to the general case, see Bj\"{o}rner and Wachs \cite[Lemma 4.2]{BW}, and eventually  finished the proof of Theorem \ref{Thm3}.

The way in which Bj\"{o}rner and Wachs \cite{BW} proved Theorem \ref{Thm3} does not seem to apply to Theorem \ref{Le1}.  To prove Theorem \ref{Le1}, we introduce  $(F,w)$-partitions of type $B$ which reduce to the $(P,w)$-partitions due to Stanley \cite{Stanley2} when $P$ is a forest  and $w$  is a signed  labeling of $P$ with labeling set $\{-1,-2,\ldots,-n\}$.  When $w$ is an ordinary labeling,  our approach leads to a direct combinatorial  proof of Theorem \ref{Thm3} of Bj\"{o}rner and Wachs.

Let $\mathbb{N}$ be the set of nonnegative integers. A  $(F,w)$-partition of type $B$ is a map $f\colon V(F) \longrightarrow
\mathbb{N}$ satisfying the following conditions

\begin{itemize}
\item[(1)] $f(x)\leq f(y)$ if $x\geq_F y$;
\item[(2)] $f(x)<f(y)$ if $x>_Fy$ and $w(x)<w(y)$;
\item[(3)] $f(u)\geq 1$ if $u$ is a root of $F$ with $w(u)>0$.
\end{itemize}
We denote by $A_B(F,w)$ the set of  $(F,w)$-partitions of type $B$. For $f\in A_B(F,w)$, let
\[ |f|=\sum_{x\in F}f(x).\]

We shall compute the generating function for the $q$-counting of
$(F,w)$-partitions of  type $B$ in two different ways.  These two different countings lead to the following identity
\begin{equation}\label{relation}
\frac{\sum_{\sigma\in\mathcal{L}(F,w)}q^{\mathrm{maj}_B(\sigma)}}{(1-q)(1-q^2)\cdots
(1-q^n)}=\frac{q^{\mathrm{maj}_B(F,w)}}{\prod_{u\in
F}(1-q^{h_u})},\end{equation} which is equivalent to the formula (\ref{Le1equ}).

\begin{lemma}\label{Lem3}
Let $F$ be a forest and $w$ a signed labeling of $F$. Then
\begin{equation}\label{hook}
\sum_{f\in A_B(F,w)}q^{|f|}=\frac{q^{\mathrm{maj}_B(F,w)}}{\prod_{u\in
F}(1-q^{h_u})}.\end{equation}
\end{lemma}

\pf For each $f\in A_B(F,w)$ and $u\in \mathrm{Des}_B(F,w)$, define $f_u\colon\ V(F) \longrightarrow
\mathbb{N}$ as
\begin{equation*}
f_u(x)=\left\{\begin{array}{ll}
f(x)-1, &\mbox{if $x\leq_F u$},\\[3pt]
f(x),&\mbox{otherwise}.\end{array}\right.
\end{equation*}
Assume $\mathrm{Des}_B(F,w)=\{u_1,u_2,\ldots,u_k\}$. Define recursively $f_{u_1u_2\cdots u_k}=(f_{u_1u_2\cdots u_{k-1}})_{u_k}$. Since $f_{uv}=f_{vu}$ then $f_{u_1u_2\cdots u_k}$ is independent of the order of $u_1,\ldots,u_k$. Thus $f_{u_1u_2\cdots u_k}$ is well defined. It is easy to see the resulting $f_{u_1u_2\cdots u_k}$ are maps from $V(F)$ to $\mathbb{N}$ satisfying only condition (1). Such maps are formally called  $P$-partitions where $P$ is a forest, see Stanley \cite[Chapter 4]{Stanley1}. In such a way,  for fixed $(F,w)$ we
establish a bijection  $\varphi\colon f\longmapsto f_{u_1u_2\cdots u_k}$ between the
following two sets\begin{equation}\label{onetoone} \varphi\colon\
A_B(F,w)\longrightarrow \{\mbox{$F$-partitions}\},\end{equation}which,
for each $f\in A_B(F,w)$, satisfies
\begin{align}|f|&=\sum_{u\in
\mathrm{Des}_B(F,\,w)}h_u+|\varphi(f)|\nonumber\\[5pt]
&=\mathrm{maj}_B(F,w)+|\varphi(f)|.\label{rela1}
\end{align}
Stanley \cite[Proposition 22.1]{Stanley2} has proved (see Sagan \cite{Sagan} for a combinatorial proof)
\begin{equation}\label{gene1}\sum_{g}q^{|g|}=\frac{1}{\prod_{u\in F}(1-q^{h_u})},\end{equation} where
the sum ranges over all $F$-partitions.
(\ref{rela1})  together with  (\ref{gene1})  implies (\ref{hook}). This completes the proof. \qed

Lemma \ref{Lem3} gives one formulation for the generating function of $A_B(F,w)$. To give the other, we go on by extending  a fundamental result concerning the usual $(P,w)$-partitions. For $\sigma=\sigma_1\cdots\sigma_n\in B_n$, a map
$f\colon \{\sigma_1,\ldots,\sigma_n\}\longrightarrow \mathbb{N}$ is said to be $\sigma$-compatible if it satisfies the following conditions
\begin{itemize}
\item[(i)] $f(\sigma_1)\geq f(\sigma_2)\geq\cdots\geq f(\sigma_n)$;
\item[(ii)] For $i\in [1,n-1]$, $f(\sigma_i)>f(\sigma_{i+1})$ if $\sigma_i>\sigma_{i+1}$;
\item[(iii)] $f(\sigma_n)\geq 1$ if $\sigma_n>0$.
\end{itemize}
Let $A_{\sigma}^B$ denote the set of all $\sigma$-compatible maps.  Then we have the following decomposition.

\begin{lemma}\label{decompsition3}
Let $F$ be a forest with a signed labeling $w$. Then
\begin{equation}\label{dec1}
A_B(F,w)=\biguplus_{\sigma\in
\mathcal{L}(F,w)}A_\sigma^B.\end{equation}
\end{lemma}

Note that if $w$ is  a signed  labeling with the labeling set $\{-1,-2,\ldots,-n\}$ then Lemma \ref{decompsition3} coincides with
 the decomposition of Stanley \cite[Theorem 6.2]{Stanley2} where the poset is a forest.

\pf  By the definition of $(F,w)$-partitions of type $B$ and the decomposition of Stanley \cite[Theorem 6.2]{Stanley2}, it is not hard to see  for each $f\in A_B(F,w)$ there exists a unique
linear extension $\sigma\in \mathcal{L}(F,w)$ such that $f\in A_\sigma^B$.
It suffices to verify that for each  $\sigma\in \mathcal{L}(F,w)$, if   $f\in A_\sigma^B$ then $f\in A_B(F,w)$. In other words, it is necessary to show that $f(\sigma_i)\geq 1$ whenever $\sigma_i>0$. Note that condition (iii) ensures this holds  in the case of $\sigma_n>0$. So we are left with  the case  $\sigma_n<0$. Assuming $\sigma_n<0$, let $\sigma_j$ be the rightmost entry of $\sigma$ such that $\sigma_j>0$. Then we have $\sigma_j>\sigma_{j+1}$ since $\sigma_{j+1}<0$. Given condition (ii), we see that $f(\sigma_j)\geq 1$. Therefore we have  $f(\sigma_i)\geq 1$ whenever $\sigma_i>0$. This completes the proof. \qed

\begin{lemma}\label{Lem2}
Let $F$ be a forest and $w$ a  signed labeling of $F$. Then
\begin{equation}\label{abf}
\sum_{f\in A_B(F,w)}q^{|f|}=\frac{\sum_{\sigma\in\mathcal{L}(F,w)}q^{\mathrm{maj}_B(\sigma)}}{(1-q)(1-q^2)\cdots (1-q^n)}.
\end{equation}
\end{lemma}

\pf For $f\in A_\sigma^B$ and
 $1\leq i\leq n-1$, let
\[ p_i=f(\sigma_i)-f(\sigma_{i+1}),\]
 and let $p_n=f(\sigma_n)$. Clearly,
\[|f|=p_1+2p_2+\cdots+np_n.\]
Under the conditions (i), (ii), and (iii), we see that $p_i>0$ if $i\in
\mathrm{Des}_B(\sigma)$ and $p_i\geq 0$ otherwise. Hence
\begin{align*}
\sum_{f\in
A_\sigma^B}q^{|f|}&=\sum_{(p_1,\ldots,p_n)}q^{p_1+2p_2+\cdots+np_n}=q^{\mathrm{maj}_B(\sigma)}\sum_{(q_1,\ldots,q_n)\atop
{q_i\geq
0}}q^{q_1+2q_2+\cdots+nq_n}\\
&=\frac{q^{\mathrm{maj}_B(\sigma)}}{(1-q)(1-q^2)\cdots(1-q^n)}.
\end{align*}
Now, using the decomposition  (\ref{dec1}), we get (\ref{abf}). \qed

Combining
Lemma
\ref{Lem3} and  Lemma \ref{Lem2}, we deduce  Theorem \ref{Le1}.
We now need one more lemma for the proof of Theorem \ref{Thm1}, which shows that the R-major index for singed permutations is equidistributed with Reiner's major index  (\ref{Rmaj}).

\begin{lemma}\label{lastlemma} For $n\geq 0$, we have
\begin{equation} \label{sbn}
\sum_{\sigma\in B_n}t^{\mathrm{p}(\sigma)}q^{\mathrm{maj}_B(\sigma)}=(1+tq)^n[n]!.
\end{equation}
\end{lemma}
\pf  We will give a bijection $\psi\colon \sigma\longmapsto \tau$ on $B_n$  such that $\mathrm{maj}_B(\sigma)=\mathrm{maj}_R(\tau)$ and $\mathrm{p}(\sigma)=\mathrm{n}_1(\tau)$. Then the formula (\ref{sbn})
 follows from Reiner's formula (\ref{rg}).  For $\sigma=\sigma_1\cdots \sigma_n\in B_n$, let $\{\sigma_{i_1},\ldots,\sigma_{i_k}\}_<$ (resp., $\{\sigma_{j_1},\ldots,\sigma_{j_{l}}\}_<$) be the set of positive (resp., negative) entries of $\sigma$, where the subscript $_<$ means $\sigma_{i_1}<\cdots<\sigma_{i_k}$. Define $\tau_{i_s}=-\sigma_{i_{k+1-s}}$ for $1\leq s\leq k$, and
$\tau_{j_t}=-\sigma_{j_{l+1-t}}$ for $1\leq t\leq l$. Let $\tau=\tau_1\cdots\tau_n$. It is clear that $\mathrm{p}(\sigma)=\mathrm{n}_1(\tau)$. Then it can be checked that  $\mathrm{Des}_B(\sigma)=\mathrm{Des}_R(\tau)$. The details are omitted.
 Thus $\psi$ is  the required bijection. This completes the proof.
\qed

We are now ready to prove Theorem \ref{Thm1}. We use a similar technique as given by Bj\"{o}rner and Wachs \cite{BW} for deriving the generating function for the major index of a given forest with the ordinary labelings.

\noindent\textit{Proof of Theorem \ref{Thm1}.} We aim to
establish the relation (\ref{main}) by evaluating the double sum
\[\sum_{(F,w)\in B_n(F)}\sum_{\sigma\in \mathcal{L}(F,w)}t^{\mathrm{p}(F,w)}q^{\mathrm{maj}_B(\sigma)}\] in two
different ways.
By Lemma \ref{Le1}, we have
\begin{align*}
\sum_{(F,w)\in B_n(F)}\sum_{\sigma\in
\mathcal{L}(F,w)}t^{\mathrm{p}(F,w)}q^{\mathrm{maj}_B(\sigma)}&=\sum_{(F,w)\in B_n(F)}t^{\mathrm{p}(F,w)}q^{\mathrm{maj}_B(F,w)}\frac{[n]!}{\prod_{u\in
P}[h_u]}\\&=\frac{[n]!}{\prod_{u\in
P}[h_u]}\sum_{(F,w)\in B_n(F)}t^{\mathrm{p}(F,w)}q^{\mathrm{maj}_B(F,w)}.
\end{align*}

On the other hand, we may compute the above double sum by exchanging the order of summation. Let $\chi$ denote the indicator function which equals $1$
when the statement is true and $0$ otherwise. Then we have
\begin{align*}
\sum_{(F,w)\in B_n(F)}\sum_{\sigma\in
\mathcal{L}(F,w)}t^{\mathrm{p}(F,w)}q^{\mathrm{maj}_B(\sigma)}&=
\sum_{(F,w)\in B_n(F)}\sum_{\sigma\in B_n}t^{\mathrm{p}(F,w)}q^{\mathrm{maj}_B(\sigma)}\chi(\sigma\in \mathcal{L}(F,w))\\
&=\sum_{\sigma\in
B_n}q^{\mathrm{maj}_B(\sigma)}\sum_{(F,w)\in B_n(F)}t^{\mathrm{p}(F,w)}\chi(\sigma\in
\mathcal{L}(F,w))\\
&=\sum_{\sigma\in
B_n}q^{\mathrm{maj}_B(\sigma)}t^{\mathrm{p}(\sigma)}\sum_{(F,w)\in B_n(F)}\chi(\sigma\in
\mathcal{L}(F,w)).
\end{align*}
Recall that for any permutation $\pi\in S_n$, \mdseries{Bj\"{o}rner and Wachs
\cite{BW} have shown that there are \[ \frac{n!}{\prod_{u\in F}h_u}\]
ordinary labelings $w$ such that $\pi\in \mathcal{L}(F,w)$. Clearly,
 this counting argument also applies to a signed permutation  $\sigma\in B_n$. Consequently,
\begin{align*}
\sum_{(F,w)\in B_n(F)}\sum_{\sigma\in\mathcal{L}(F,w)}t^{\mathrm{p}(F,w)}q^{\mathrm{maj}_B(\sigma)}&=\sum_{\sigma\in
B_n}q^{\mathrm{maj}_B(\sigma)}t^{\mathrm{p}(\sigma)}\frac{n!}{\prod_{u\in F}h_u}\\
&=\frac{n!}{\prod_{u\in F}h_u}\sum_{\sigma\in
B_n}q^{\mathrm{maj}_B(\sigma)}t^{\mathrm{p}(\sigma)}.
\end{align*}
By Lemma \ref{lastlemma}, we get
\[\sum_{(F,w)\in B_n(F)}\sum_{\sigma\in\mathcal{L}(F,w)}t^{\mathrm{p}(F,w)}q^{\mathrm{maj}_B(\sigma)}
=\frac{n!}{\prod_{u\in F}h_u}(1+tq)^n[n]!.\]

The above double counting gives the following relation
\[\frac{[n]!}{\prod_{u\in
F}[h_u]}\sum_{(F,w)\in B_n(F)}t^{\mathrm{p}(F,w)}q^{\mathrm{maj}_B(F,w)}=\frac{n!}{\prod_{u\in
F}h_u}(1+tq)^n[n]!,\]which is equivalent to (\ref{main}). This completes the proof.
\qed

Based on Theorem \ref{Thm1}, it is easy to derive Theorem \ref{Tnmaj}.

\noindent\textit{Proof of Theorem \ref{Tnmaj}.} Setting $q\rightarrow q^2$ and $t\rightarrow q^{-1}$ in  (\ref{main}), the left-hand side  becomes
\begin{align*}
\sum_{(F,w)\in\, B_n(F)}(q^{-1})^{\mathrm{p}(F,w)}q^{2\mathrm{maj}_B(F,w)}
&=\sum_{(F,w)\in\, B_n(F)}q^{2\mathrm{maj}_B(F,w)-\mathrm{p}(F,w)}\\[5pt]
&=\sum_{(F,w)\in\, B_n(F)}q^{\mathrm{rmaj}(F,w)}
\end{align*} and the right-hand side can be written as
\begin{align*}
\frac{n!}{\prod_{u\in
F}h_u}(1+q^{-1}q^2)^n \prod_{u\in F}[h_u]_{q^2}=\frac{n!}{\prod_{u\in
F}h_u}\prod_{u\in F}[2h_u].\end{align*}
This completes the proof.\qed

\section{A correspondence}

From Theorem \ref{Tnmaj} and Theorem \ref{Tfmaj}, one sees that the R-major index and
the flag major index are equidistributed for signed labeled forests.
One is naturally led to the question of finding a correspondence that
explains the equidistribution property. This is the objective of this section
to provide such a correspondence. Of course, this bijection can be considered
as an alternative proof of Theorem \ref{Tnmaj}.

 Define  a bijection $\phi\colon\ (F,w)\longmapsto
(F,w')$ on $B_n(F)$ as follows. For each vertex $u\in F$,
\begin{itemize}
\item[(1)] $w'(u)$ has the same sign with $-w(u)$,
\item [(2)] $|w'(u)|=n+1-|w(u)|$.
\end{itemize}

\begin{theorem}
The above map $\phi$ is a bijection with the following property \[\mathrm{rmaj}(F,w)=\mathrm{fmaj}(F,w').\]
\end{theorem}

\pf The theorem  holds if we can show that
\begin{equation}\label{relaf-p}
\mathrm{maj}_B(F,w)=\mathrm{maj}(F,w')+\mathrm{p}(F,w)\ \ \mathrm{and}\ \ \mathrm{p}(F,w)=\mathrm{n}(F,w').\end{equation}
We proceed to prove (\ref{relaf-p}) by induction on the number of vertices of $F$.
Without loss of generality, we may assume that  $F$ is a tree.

If $F$ has only one vertex,  it is easy to check (\ref{relaf-p}). So we
 may assume that $F$ has at least two vertices. Let $u_0$ be the root of $F$, and $C(u_0)$ the set of children of $u_0$. By the definition of $\phi$, we see that
  $\mathrm{p}(F,w)=\mathrm{n}(F,w')$. Now we claim that
\begin{equation}\label{ggg}
\mathrm{maj}_B(F,w)=\mathrm{maj}(F,w')+\mathrm{p}(F,w).
\end{equation}
Here are two cases.

Case 1: $w(u_0)>0$.  We partition $C(u_0)$ into the following three subsets.
\begin{align*}
P_>&=\{u\in C(u_0)\colon\, w(u)>w(u_0)\},\\[5pt]
P_<&=\{u\in C(u_0)\colon\, 0<w(u)<w(u_0)\},\\[5pt]
N&=\{u\in C(u_0)\colon\, w(u)<0\}.
\end{align*}
For any vertex $u\in F$, let $F_u$ be the subtree of $F$ rooted at $u$. Then
\begin{equation*}
\begin{split}
&\mathrm{Des}_B(F,w)\\[5pt]
=&\left(\bigcup_{u\in P_>}\mathrm{Des}_B(F_u,w_u)\right)\bigcup\left(\bigcup_{u\in P_<}\mathrm{Des}_B(F_u,w_u)\backslash P_<\right)\bigcup\left(\bigcup_{u\in N}\mathrm{Des}_B(F_u,w_u)\right)\bigcup\,\{u_0\},
\end{split}
\end{equation*}
where $w_u$ is the signed labeling of $F_u$ induced by $w$. So we have
\begin{equation}\label{majBB}
\begin{split}
\mathrm{maj}_B(F,w)=&\sum_{u\in P_>}\mathrm{maj}_B(F_u,w_u)+\sum_{u\in P_<}\mathrm{maj}_B(F_u,w_u)\\[5pt]
&+\sum_{u\in N}\mathrm{maj}_B(F_u,w_u)+n-\sum_{u\in P_<}h_u.
\end{split}
\end{equation}

Let us compute the major index of $w'$.  Let
\begin{align*}
N_>&=\{u\in C(u_0)\colon\, 0>w'(u)>w'(u_0)\},\\[5pt]
N_<&=\{u\in C(u_0)\colon\, w'(u)<w'(u_0)\},\\[5pt]
P&=\{u\in C(u_0)\colon\, w'(u)>0\}.
\end{align*}
Then
\begin{equation*}
\begin{split}
&\mathrm{Des}(F,w')\\[5pt]
=&\left(\bigcup_{u\in N_>}\mathrm{Des}_B(F_u,w'_u)\bigcup N_>\right)\bigcup\left(\bigcup_{u\in N_<}\mathrm{Des}(F_u,w'_u)\right)\bigcup\left(\bigcup_{u\in P}\mathrm{Des}(F_u,w'_u)\bigcup P\right),
\end{split}
\end{equation*}
from which we deduce that
\begin{equation}\label{gg}
\begin{split}
\mathrm{maj}(F,w')=&\sum_{u\in N_>}\mathrm{maj}(F_u,w')+\sum_{u\in N_<}\mathrm{maj}(F_u,w')\\[5pt]
&+\sum_{u\in P}\mathrm{maj}(F_u,w')+\sum_{u\in N_>}h_u+\sum_{u\in P}h_u.
\end{split}
\end{equation}
By the definition of $\phi$, it is not hard to verify
\begin{equation}\label{setset}
N_>=P_>,\ \ N_<=P_<\ \ \mathrm{and}\ \ P=N.
\end{equation}
Therefore, by  (\ref{majBB}), (\ref{gg}), (\ref{setset}), and the induction hypothesis, we find that
\allowdisplaybreaks
\begin{align*}
\mathrm{maj}_B(F,w)=&\sum_{u\in P_>}\left(\mathrm{maj}(F_u,w'_u)+\mathrm{p}(F_u,w_u)\right)+\sum_{u\in P_<}\left(\mathrm{maj}(F_u,w'_u)+\mathrm{p}(F_u,w_u)\right)\\[5pt]
&+\sum_{u\in N}\left(\mathrm{maj}(F_u,w'_u)+\mathrm{p}(F_u,w_u)\right)+n-\sum_{u\in P_<}h_u\\[5pt]
=&\sum_{u\in N_>}\mathrm{maj}(F_u,w'_u)+\sum_{u\in N_<}\mathrm{maj}(F_u,w'_u)
+\sum_{u\in P}\mathrm{maj}(F_u,w'_u)\\[5pt]
&+\mathrm{p}(F,w)-1+n-\sum_{u\in N_<}h_u.\\[5pt]
=&\mathrm{maj}(F,w')+\mathrm{p}(F,w)-1+n-\sum_{u\in N_<}h_u-\sum_{u\in N_>}h_u-\sum_{u\in P}h_u\\[5pt]
=&\mathrm{maj}(F,w')+\mathrm{p}(F,w)-1+n-\sum_{u\in C(u_0)}h_u,
\end{align*}
which reduces to (\ref{ggg}) since
\[\sum_{u\in C(u_0)}h_u=n-1.\]

Case 2:   $w(u_0)<0$. We can use a similar argument to that for Case 1.
So we reach the conclusion that  (\ref{relaf-p}) holds for any tree. This completes the proof.\qed

\section{Concluding remarks}

We conclude this paper with two questions.
 While we have derived the generating functions
 for the flag major index and the R-major index of signed labeled
forests, it would be interesting to give a suitable definition
  of the negative major index (\ref{nmaj}) for signed labeled forests.
Intuitively,  a  natural choice would be
\begin{equation}\label{snmaj1}
\mathrm{nmaj}(F,w)=\mathrm{maj}(F,w)+\mathrm{n}_1(F,w)+\mathrm{n}_2(F,w).
\end{equation}
However, the above statistic is not  equidistributed with  the inversion number (\ref{SLFinv}).

For type $D_n$ permutations, Biagioli \cite{Biagioli2} defined the  major index as
follows
\begin{equation}\label{dmaj}
\mathrm{dmaj}(\sigma)=\mathrm{maj}(\sigma)+\mathrm{n}_2(\sigma),
\end{equation} and shown it is equidistributed with the length function of $D_n$. However, the following statistic
\begin{equation}\label{snmaj2}
\mathrm{dmaj}(F,w)=\mathrm{maj}(F,w)+\mathrm{n}_2(F,w)
\end{equation}
for labeled forests of type $D_n$   is not equidistributed with
the inversion number (\ref{ESLFinv}) of $D_n(F)$.
 We would pose the question of finding an appropriate major index
 for $D_n(F)$ which is equidistributed with the inversion number.

\end{document}